\documentclass[a4paper,11pt,numreferences,mathsec,kaplist]{isueps}
\usepackage{isu}

\begin{document}
\setcounter{aqwe}{2} 
\begin{article}
\begin{opening}
\udk{51-74}

\title{Control of accuracy on Taylor-collocation method for load leveling problem\thanks{This research was supported by RFBR Grant No. 18-31-00206.}}

\author{Samad~\surname{Noeiaghdam}}
\institute{Baikal School of BRICS, Irkutsk National Research
Technical University, Irkutsk, Russian Federation\\
South Ural State University, Chelyabinsk, Russian Federation}

\author{Denis~\surname{Sidorov}}
\institute{Melentiev Energy Systems Institute, SB of RAS, Irkutsk, Russian Federation}

\author{Ildar~\surname{Muftahov}}
\institute{Melentiev Energy Systems Institute, SB of RAS, Irkutsk, Russian Federation\\
Irkutsk Division of Main Computer Centre, JSC Russian Railways, Irkutsk, Russian Federation}

\author{Aleksei~\surname{Zhukov}}
\institute{Institute of Solar-Terrestrial Physics SB RAS, Irkutsk, Russian Federation\\
Melentiev Energy Systems Institute, SB of RAS, Irkutsk, Russian Federation}

\runningtitle{TAYLOR-COLLOCATION METHOD FOR LOAD LEVELING PROBLEM}
\runningauthor{S. Noeiaghdam, D. Sidorov, I. Muftahov, A. Zhukov}

\begin{abstract}
High penetration of renewable energy sources coupled with
decentralization of transport and heating loads in future power
systems will result even more complex unit commitment problem
solution using energy storage system scheduling for efficient load
leveling. This paper employees an adaptive approach to load leveling
problem using the Volterra integral dynamical models. The problem is
formulated as solution of the Volterra integral equation of the
first kind which is attacked using Taylor-collocation numerical
method which has the second-order accuracy and enjoys
self-regularization properties, which is associated with confidence
levels of system demand. Also the CESTAC method is applied to find
the optimal approximation, optimal error and optimal step of
collocation method. This adaptive approach is suitable for energy
storage optimization in real time. The efficiency of the proposed
methodology is demonstrated on the Single Electricity Market of the
Island of Ireland.
\end{abstract}


\end{opening}

\avtogle{S. Noeiaghdam, D. Sidorov, I. Muftahov, A. Zhukov}{Control of accuracy on Taylor-collocation method for load leveling problem}

\section{Introduction}
Robust numerical methods design for optimized energy storage and
load leveling is one of the most imminent challenges in power
systems. The conventional least-cost dispatch of available
generation to meet the forecasted load will no longer be suitable
for the purposes because of observed decentralization of power
systems including electrical transport and heating loads and
renewable energy use worldwide. Various methods have been examined
to solve this challenging problem including the evolutionary
methods, and genetic algorithms, Lagrangian relaxation, mixed
integer linear programming and particle swarm optimization.
The transmission system operator in Republic of
Ireland restricts the instantaneous proportion of total generation
allowed from wind wind turbines, to 50 \% maximum, in order to
maintain sufficient system inertia \cite{001}. This may result in
wind curtailment at any time.
The reduction of the uncertainty associated
with wind power can be achieved using the state of the art
forecasting methods based on contemporary machine learning theory.
Wind power and industrial and residential electric load forecasting,
over lead times of up to 48 hours, are critical for the market
operator. These forecasts are used to create day-ahead unit
commitment and economic dispatch schedules. In fact, many
transmission system operators in Russia and in other countries also
employ the shorter-term wind forecasts to draw upon system reserves
for short-term balancing. It is to be noted here that a reduction in
the forecast errors by a fraction of a percent will lead to a
substantial increase in trading profits. The efficient forecast
algorithms increase the value of wind generation. In \cite{002} it
is shown that an increase of only one percent in forecast error
caused an increase of 10 million GBP in operating costs per year for
one electric utility in the UK.
For
the detailed review of China coal-fired power units peak regulation
with a detailed presentation of the installed capacity, load
leveling operation mode readers may refer to \cite{003}. The
effectiveness of wind turbine, energy storage, and demand response
programs in the deterministic and stochastic circumstances and
influence of uncertainties of the wind, price, and demand are
assessed in the Energy Hub planning in \cite{004}. The unit
commitment solver is supposed to serve as an integral part of Energy
Hub for control of the interconnections of heterogeneous energy
infrastructures, including non-synchronous renewable sources. This
paper presents a methodology for load leveling using dynamical
models based on efficient numerical solution of the integral
equations (IE) of the first kind. We consider the future power
systems with battery energy storages of various efficiencies. The
load forecasts over lead times of 24 hours is used. The load
leveling problem is formulated as inverse problem. These integral
equations efficiently solve such an inverse problem given both the
time-dependent efficiencies and the availability of
generation/storage technologies. Electric load can be forecasted
using various mathematical models including classical statistical
methods of time series analysis, regressive models or advanced
machine learning methods including deep learning. The evolutionary
character of power systems when future values depend upon previous
values can be described by evolutionary integral equations which are
also known as Volterra equations. The integrand (or kernel) of
employed equations has 1st kind discontinuities along the continuous
curves starting at the origin. Such piecewise continuous kernels
K(t,s) take into account both efficiencies of the different storing
technologies and their proportions (which could be time-dependent)
in the total power generation. Efficiencies of the different
storages may depend not only on their age and use duration but also
on the state of charge, see \cite{005}. This phenomena can be
modeled using the nonlinear Volterra integral equations.

In recent years many numerical or semi-analytical methods have been
presented to solve the linear and non-linear first kind IEs
\cite{man2,
man1,
man7,siz5,siz3,siz1}. In this paper, we
develop further the numerical methods for solution of such integral
equations and continue our previous work \cite{006} where we
demonstrated how these novel models can be applied for storage
modeling. In this work, the Taylor-collocation method
\cite{16,15,18,
17
} is considered to solve the mentioned
problem based on the stochastic arithmetic and the approximate
results are validated by using the CESTAC method. Also, the optimal
iteration and the optimal approximation of method are found.
Recently, the CESTAC method has been applied to implement the
numerical methods for finding the approximate solution of different
problems \cite{008,009,
0010,0011,007}. In this method, instead
of using the mathematical packages such as Matlab, Mathematica and
the others, the CADNA library is applied. Also, in this library the
logical programs can be written by statements of C/C++, FORTRAN or
ADA \cite{cad2,cad6}. Some of the advantages of using the
CESTAC method and CADNA library are:

\begin{itemize}
	\item In the CESTAC method, not only the optimal numerical solution can
be produced but also, the optimal iteration can be obtained
\cite{008,009,0010,007}.

	\item The CADNA library is able to detect any instability in
mathematical operations, branching, functions and so on but the FPA
has not these abilities. 

	\item In the FPA, the termination criterion depends on a small
parameter like $h$ or epsilon. For epsilon enough large,  the
iterations can be stopped before finding the suitable approximation
and for small values of epsilon the unnecessary iterations can be
produced without improving the accuracy of the results
\cite{cad9,cad12}. In the SA, the numerical results do not depend on
the value epsilon and a new stopping condition is replaced which is
independent of epsilon and existence of exact solution.
\end{itemize}

This paper is structured into four sections. The background on
Volterra integral equations theory and numerical methods (including
results on convergence) is described in Section 2. Section 3
describes the results of the numerical experiments on real data from
the Irish test system. Section 4 draws a brief conclusion

\section{Volterra model}

The proposed collocation-type numerical method has the second-order
accuracy and enjoys self-regularization properties as mesh step is
associated with confidence levels of system demand. The proposed
approach is suitable for energy storage control in real time.

Let us consider the  first kind classical Volterra IE
\begin{equation}\label{1}
\int_{0}^{t} K(t,s) x(s) ds  = f(t),~ 0 \leq s \leq t \leq T,
\end{equation}
where $K(t,s)$ is a discontinuous kernel along continuous curves
$\alpha_i, i=1,2,\cdots,m-1$ as
\begin{equation}\label{2}\left\{
\begin{array}{l}
  K_1(t,s),~~~~~ 0= \alpha_0(t)  \leq s \leq  \alpha_1(t),\\
    \\
  K_2(t,s),~~~~~ \alpha_1(t)  \leq s \leq  \alpha_2(t),\\
    ~~~\vdots\\
    K_m(t,s),~~~~~ \alpha_{m-1}(t)  \leq s \leq  \alpha_m(t) =r \leq 1. \\
\end{array}\right.
\end{equation}
Finally,  Eq. (\ref{1}) can be written in the following form
$$
\int_{\alpha_0(t)}^{\alpha_1(t)} K_1(t,s) x(s) ds +
\int_{\alpha_1(t)}^{\alpha_2(t)} K_2(t,s) x(s) ds + \cdots +
$$
\begin{equation}\label{3}
+\int_{\alpha_{m-1}(t)}^{\alpha_m(t)} K_m(t,s) x(s) ds = f(t),
\end{equation}
where $f(0)=0$. Also, the compacted scheme of Eq. (\ref{3}) is given
as
\begin{equation}\label{4}
\sum_{p=1}^m \int_{\alpha_{p-1}(t)}^{\alpha_p(t)} K_p(t,s) x(s) ds =
f(t).
\end{equation}

Recently, the IE with jump discontinuous kernel (\ref{4}) have been
solved by many mathematical methods
\cite{dav,0015,dn8,dn5,
006,0013,0014,
0025} based on the
floating point arithmetic. In these researches and many other
articles, the validation of presented methods are investigated by
using the absolute error. But does this tool is a proper instrument
to validate the numerical results? In many researches we can find
the dependence of obtained results to the exact solution as
\begin{equation}\label{cond1}
\left| x(s) - x_n(s) \right| \leq \varepsilon.
\end{equation}
It means that without having the exact solution we can not study the
accuracy of methods. Also, how are we sure about optimal value for
number of iteration $n$? On the other hand, how can we apply the
condition (\ref{cond1}) without knowing the proper value of
$\varepsilon$.

In this study, in order to solve the mentioned problems which are
based on the FPA, the Taylor-collocation method based on the
stochastic arithmetic is applied and then the obtained results are
validated by using the CESTAC method and the CADNA library.
Moreover, a novel stopping condition
\begin{equation}\label{cond2}
\left| x_{n+1} (s) -  x_n (s) \right| = @.0,
\end{equation}
is presented instead of criterion (\ref{cond1}) that it depends on
two successive approximations. New algorithm will be stopped when
difference of two successive approximations equals to informatical
zero $@.0$. Let
\begin{equation}\label{5}
 x_n(s) = \sum_{j=0}^n \frac{1}{j!} x^{(j)} (c) (s-c)^j + O(h^{n+1}),
\end{equation}
be the $n$-th order of Taylor polynomial at point $s=c$ where
\begin{equation}\label{2222}
e_n(s) = |x(s)-x_n(s)| = O(h^{n+1}).
\end{equation}
By putting Eq. (\ref{5}) into Eq. (\ref{4}) we get
\begin{equation}\label{6}
\sum_{p=1}^m \int_{\alpha_{p-1}(t)}^{\alpha_p(t)} K_p(t,s) \sum_{j=0}^n
\frac{1}{j!} x^{(j)} (c) (s-c)^j ds = f(t).
\end{equation}
Now, the collocation points
\begin{equation}\label{2-1}
r_i=a+(\frac{b-a}{n})i, ~~~~~i=0,1,\cdots, n,
\end{equation}
should be substituted in Eq. (\ref{6}) as follows
\begin{equation}\label{2-2}
\sum_{j=0}^n \frac{1}{j!} \left[ \sum_{p=1}^m
\int_{\alpha_{p-1}(r_i)}^{\alpha_p(r_i)} k_p(r_i,s)
  (s-c)^j ds \right] x^{(j)} (c) = f(r_i).
\end{equation}

The matrix form of Eq. (\ref{2-2}) is presented as
\begin{equation}\label{2-3}
AV=F,
\end{equation}
where
$$
A=\left[
\begin{array}{cccc}
  A_{00} & A_{01} & \cdots & A_{0n} \\
  A_{10} & A_{11} & \cdots & A_{1n} \\
  \vdots & \vdots & \ddots & \vdots \\
  A_{n0} & A_{n1} & \cdots & A_{nn}
\end{array}
\right]_{(n+1)(n+1)},
$$
that
$$
A_{ij} = \sum_{p=1}^m \int_{\alpha_{p-1}(r_i)}^{\alpha_p(r_i)}
k_p(r_i,s)
  (s-c)^j ds,
$$
and
$$
V=\left[ \begin{array}{cccc}
           \bar{x}^{(0)}(c) & \bar{x}^{(1)}(c) & \cdots & \bar{x}^{(n)}(c)
         \end{array}
   \right]^T,
$$
$$
F=\left[ \begin{array}{cccc}
           F(r_0) & F(r_1)  & \cdots & F(r_n)
         \end{array}
   \right]^T.
$$

By solving system (\ref{2-3}), the coefficients $\bar{x}^{(j)}(c)$
can be found uniquely. Thus the unique solution of Eq. (\ref{4}) can
be calculated by
\begin{equation}\label{8}
 \bar{x}_n(s) = \sum_{j=0}^n \frac{1}{j!} \bar{x}^{(j)} (c) (s-c)^j.
\end{equation}



\subsection{The CESTAC methodology and the CADNA library}


Let $F$ be the set of representable values that produced by
computer. Then we can generate $G \in F$ with $P$ mantissa bits of
the binary FPA for $g \in \mathbb{R}$ as
\begin{equation}\label{0019}
G=g-\chi 2^{E-P}\gamma,
\end{equation}
where $\chi$ is the sign of $g$, $2^{-P}\gamma$ is the missing
segment of the mantissa because of round-off error and $E$ shows the
binary exponent of the result. Also we can choose two values $P =
24, 53$ for single and double accuracies
\cite{cad2,cad6}. 

Let $\gamma$ be the stochastic variable that uniformly distributed
on $[-1, 1]$. Now the perturbation on final mantissa bit of $g$ can
be made. So the casual variables with mean $(\mu)$ and the standard
deviation $(\sigma)$ can be obtained for results of $G$. We should
note that parameters $(\mu)$ and $(\sigma)$ have main rule to
determine the precision of $G$ \cite{cad9,
cad12}. By $l$ times
doing the mentioned process for $G_i , i = 1, . . . , l$ we can make
the quasi Gaussian distribution for them where mean of these values
equals to the exact $g$. Also values $\mu$ and $\sigma$ can be found
based on these $l$ samples. Algorithm 1, is introduced based on the
CESTAC method where $\tau_{\delta}$ is the value of $T$ distribution
with $l-1$ degree of freedom and confidence interval $1-\delta$.

\textbf{Algorithm 1:}

\texttt{Step 1- Provide $l$ samples of $G$ as $ \Phi = \left\{ G_1,
G_2, ..., G_l \right\}$ by means of the perturbation of the last bit
of mantissa.}

\texttt{Step 2- Find $\displaystyle G_{ave} = \frac{\sum_{i=1}^l
G_i}{l}$.}

\texttt{Step 3- Compute $\displaystyle \sigma^2= \frac{\sum_{i=1}^l
( G_i - G_{ave})^2}{l-1}$. }

\texttt{Step 4- Apply $\displaystyle C_{G_{ave}, G} =
\log_{10}\frac{\sqrt{l} \left|G_{ave} \right|}{\tau_{\delta}
\sigma}$ to show the NCSDs\\between $G$ and $G_{ave}$.}

\texttt{Step 5- If $C_{G_{ave}, G} \leq 0$ or $G_{ave}=0,$
then write $G=@.0$.}\\

Applying the mathematical methods based on the SA has many
advantages that using the mathematical methods based on the FPA. In
the SA, instead of applying some packages like Mathematica, Maple
and many others we use the CADNA library which is a novel software.
Unlike other softwares we should run it on LINUX operating system
and the CADNA commands are based on the C, C++, FORTRAN or ADA codes
\cite{008,009,
0011,007}.

By using the CESTAC method we can apply the new stopping condition
(\ref{cond2}) instead of (\ref{cond1}) which depends on two
successive approximations. Also the new condition independence of
exact solution and the tolerance value $\varepsilon$. Furthermore,
new condition will be stopped when the NCSDs of numerical results
equals to informatical zero $@.0$.

Finally, the main capability of the CESTAC method is to find the
optimal factors of numerical methods like optimal approximation,
optimal step  and optimal error of method. Also, the other main
ability of the CADNA library is to detect the numerical
instabilities \cite{cad2,cad6
}. The general
sample CADNA library program is presented in the following form\\
\texttt{$\sharp$include <cadna.h>}\\
\texttt{cadna$_{-}$init(-1);\\} \texttt{ main()\\}
\texttt{\{} \\
   \texttt{double$_{-}$st VALUE;} \\
\texttt{do}\\
    \texttt{\{}\\
The Main Program;\\
\texttt{printf(" \%s  ",strp(VALUE));}\\
\texttt{\}}\\
\texttt{ while(x[n]-x[n-1]!=0);\\}
 \texttt{cadna$_{-}$end();\\}
 \texttt{\}}


Now, we need to prove a theorem to demonstrate the equality  of the
NCSDs between $x_{n}(t), x(t)$ and $x_{n}(t),x_{n+1}(t)$.

\begin{definition} \label{def2} \cite{cad2,cad6}
For numbers $z_1, z_2 \in \mathbb{R}$, the NCSDs can be
obtained as follows\\
 (1) for $z_1\neq z_2$,
\begin{equation}\label{13}
C_{z_1,z_2}=\log_{10}\left|\frac{z_1+z_2}{2(z_1-z_2)} \right| =
\log_{10}\left|\frac{z_1}{z_1-z_2} - \frac{1}{2} \right|,
\end{equation}
(2) for all real numbers $z_1$, $C_{z_1,z_1} = +\infty$.
\end{definition}

\begin{theorem}\cite{007} \label{th6}
Assume that $x(t)$ and $x_n(t)$ are the exact and numerical
solutions of problem (\ref{4}) then
\begin{equation}\label{14}
\displaystyle C_{x_{n}(t),x(t)} - C_{x_{n}(t),x_{n+1}(t)} =
\mathcal{O}\left( h^{n+1}\right),
\end{equation}
where $C_{x_{n}(t),x(t)}$ shows the NCSDs of $x_{n}(t), x(t)$ and
$C_{x_{n}(t),x_{n+1}(t)}$ is the NCSDs of two iterations
$x_{n}(t),x_{n+1}(t)$.
\end{theorem}


\section{Results of numerical experiments on real data}

\begin{figure}[htbp]
	\centering
	\includegraphics[width=0.7\linewidth]{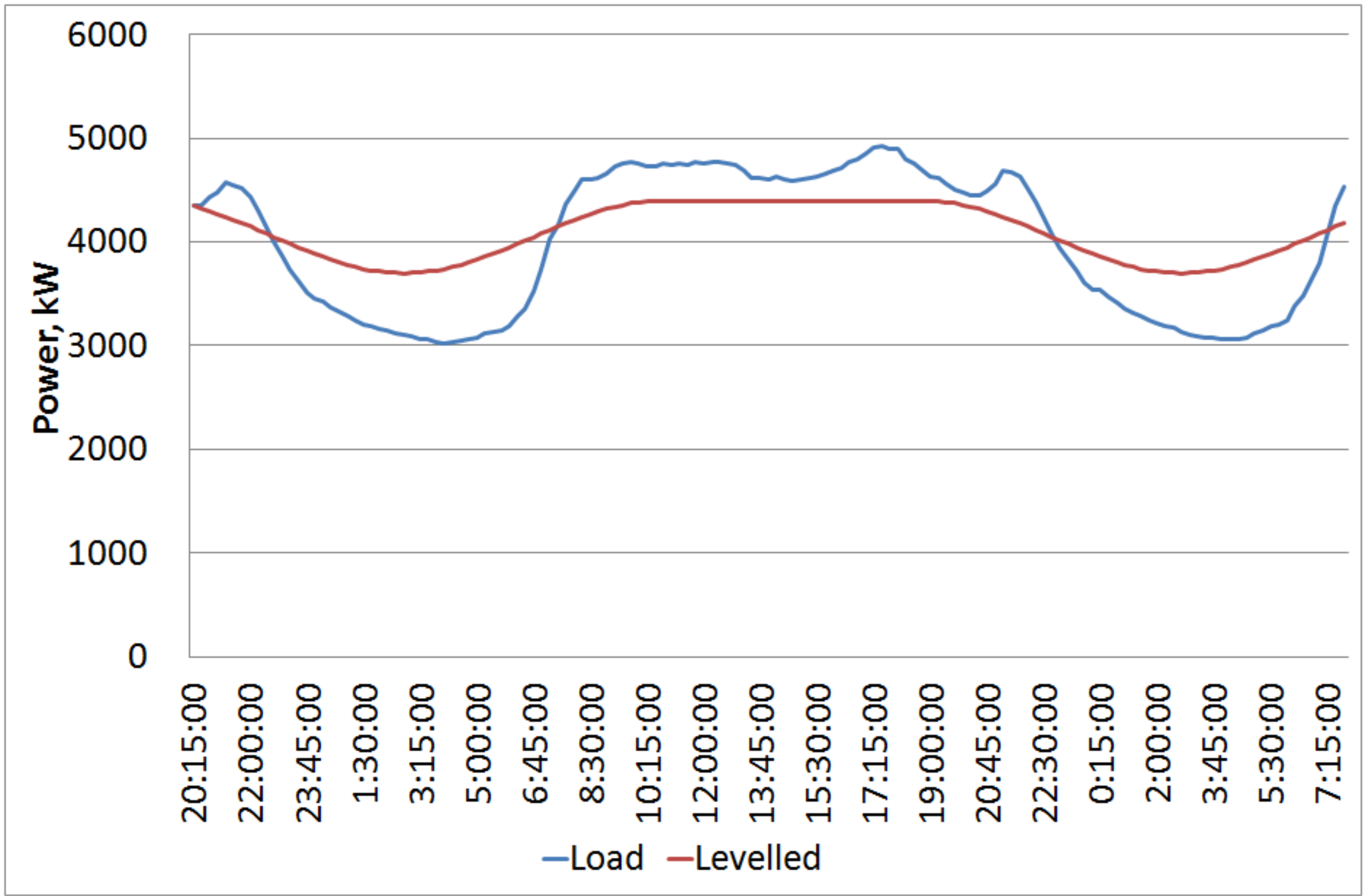}
	\caption{Leveled load on Ireland data}
	\label{fig:load}
\end{figure}

As shown in \cite{0025}, Volterra models can be used to calculate a charge and discharge strategy for storage devices.
In this paper, we use the same data (figure \ref{fig:load}) from Ireland \cite{001} to compare the Taylor-collocation method and spline collocations in the load leveling problem.

In the calculations of the charging/discharging strategy, we used a piece-wise specified core for the efficiency of the storage components:
\begin{equation*}
    K(t, s) = \left\{ \begin{array}{ll}
         \mbox{$1, \,\, 0<s<t/4,$} \\
         \mbox{$0.9, \,\, t/4<s<3 t/4,$} \\
         \mbox{$0.85, \,\, 3t/4<s<3 t.$} \\
        \end{array} \right.
\end{equation*}

\begin{figure}[htbp]
	\centering
	\includegraphics[width=0.7\linewidth]{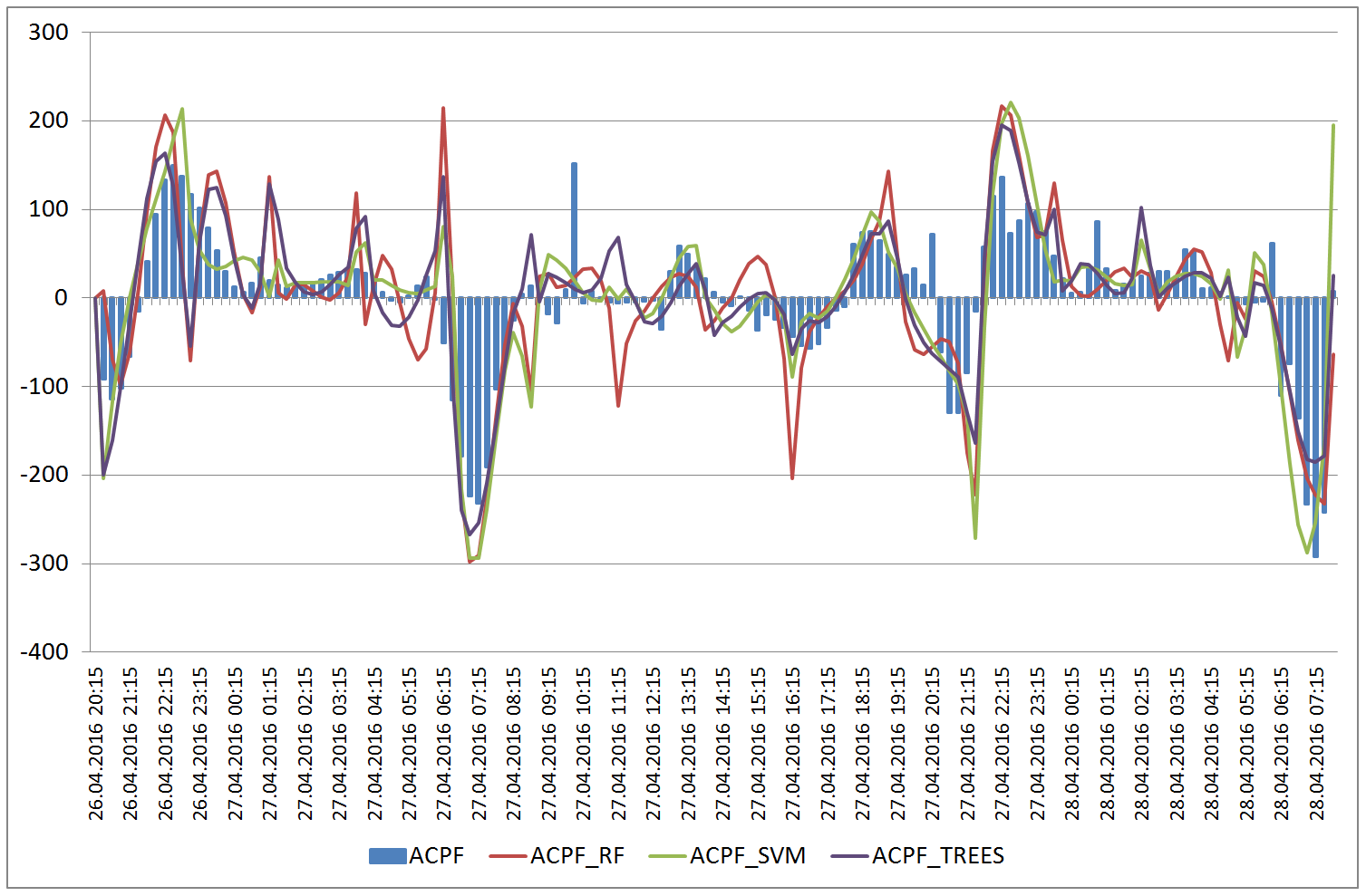}
	\caption{Charge/discharge strategies using various 24- hour forecasts and calculated by Taylor-collocation method}
	\label{fig:acpf}
\end{figure}

Since the data on the load and its forecasts are given in tabular form, we used the approximations of load and its forecasts by a 4th degree polynomials by 10 points for solving Volterra integral equation by Taylor collocations.
The result of the work on accurate and forecast data is shown in the figure (\ref{fig:acpf}).
Here alternating changing power functions (ACPF) show charge and discharge strategies.

\begin{figure}[htbp]
	\centering
	\includegraphics[width=0.7\linewidth]{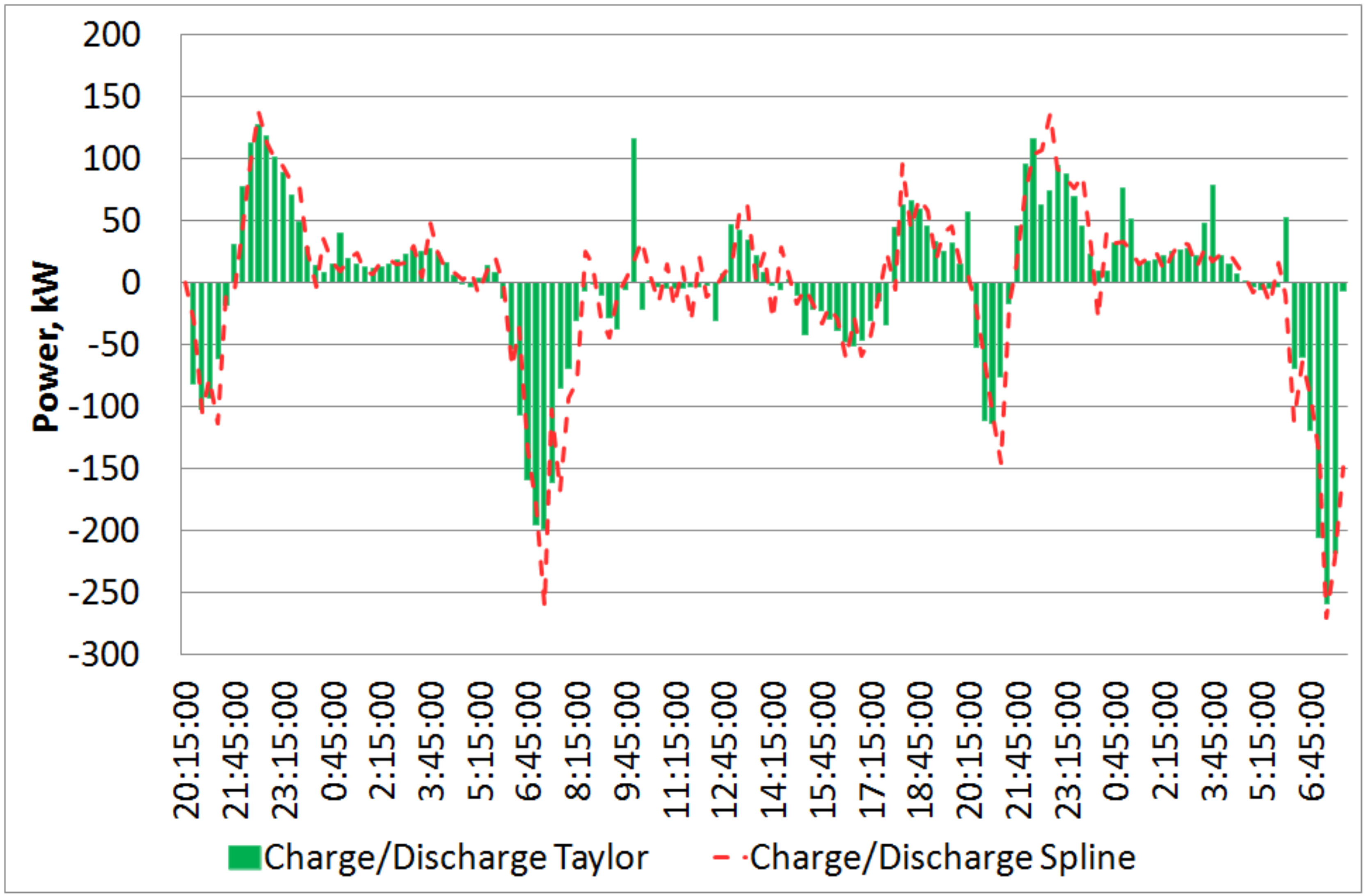}
	\caption{Charge/discharge strategies based on Taylor-collocation method and strategy  spline collocations (bars)}
	\label{fig:collocations}
\end{figure}

The Taylor collocation method has similar results to the spline collocation method as shown in the figure (\ref{fig:collocations}).

\begin{figure}[htbp]
	\centering
	\includegraphics[width=0.7\linewidth]{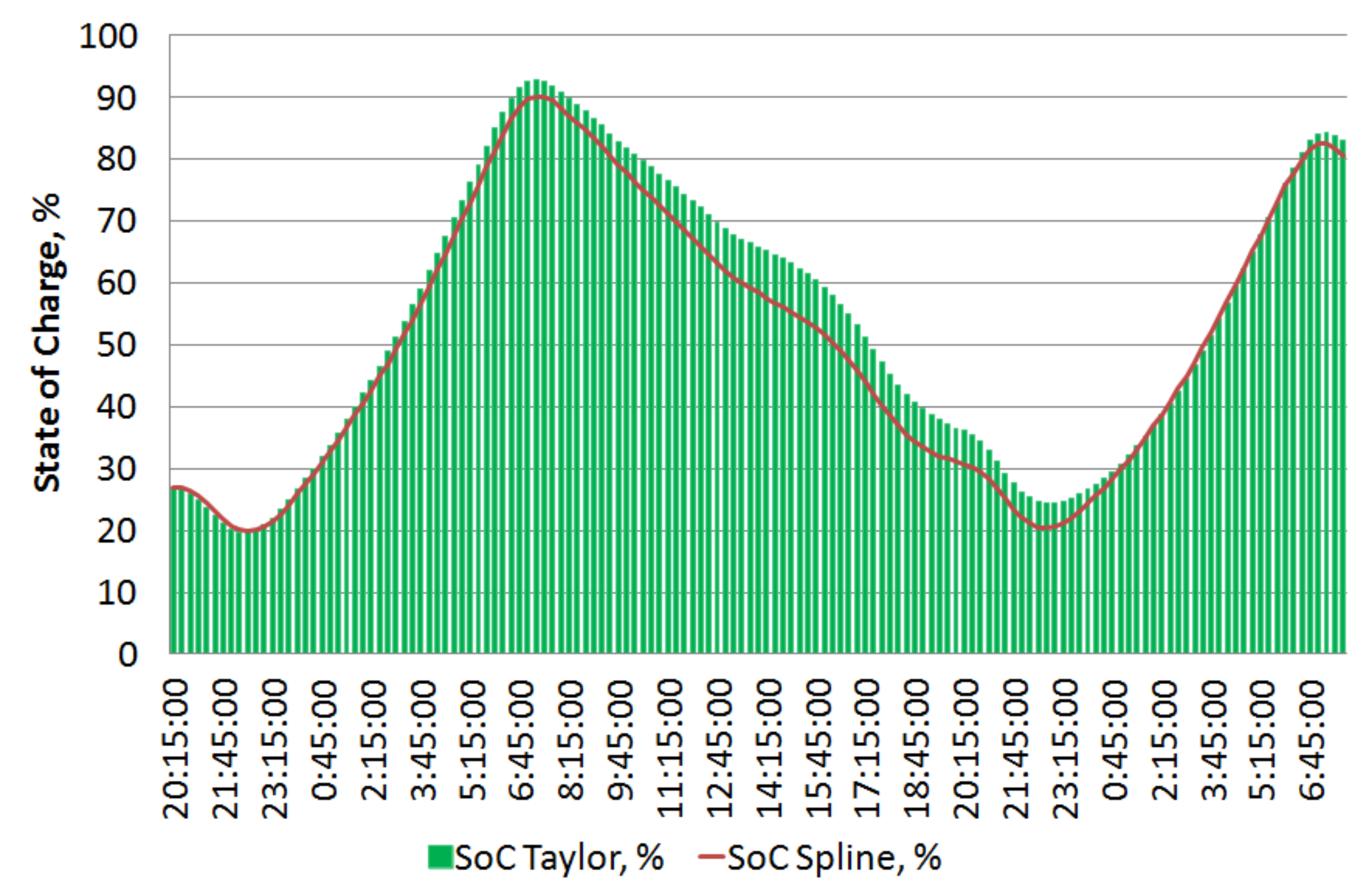}
	\caption{State of charge calculated by Taylor-collocation method and spline collocations}
	\label{fig:SoC}
\end{figure}

The state of charge (SoC) calculated by these methods is shown in the figure (\ref{fig:SoC}).
The difference between SoC values can be explained by the polynomial approximation of the right side of the Volterra integral equation for the Taylor-collocation method.

\selectlanguage{english}

\begin{table}
\caption{ Numerical results based on the CESTAC method for $f_8$.
}\label{tb20}\centering
\begin{tabular}{|l|l|l|}
  \hline
$n$&$v_n(t)$&$|v_n-v_{n+1}(t)|$  \\
    \hline
2&0.10082463733333E+003&0.10082463733333E+003\\
3&0.10841798623809E+003&0.7593348904761E+001\\
4&0.1122562709626E+003&0.38382847245E+001\\
5&0.11218760890E+003&0.68662057E-001\\
6&0.1121987347E+003&0.111258E-001\\
7&0.112198015E+003&0.7190E-003\\
8&0.1121979E+003&0.9E-004\\
9&0.112197E+003&@.0\\
  \hline
  \end{tabular}
\end{table}

In order to apply the CESTAC method to find the optimal
approximation, the optimal error and the step of Taylor-collocation
method, one of inputed values $f_i$ is chosen. Then the presented
method is applied to solve the Eq. (\ref{4}) based on the CESTAC
method. In this example, the 8-th value of $f_i$ is chosen and the
numerical results of Table \ref{tb20} are obtained. According to
these results, the optimal step of charging/discharging strategies
to find the 8-th value of them is $n=9$, the optimal approximation
is $v_9 = 0.112197E+003$ and the optimal error is $0.9E-004$.
According to Theorem 2, the NCSDs of exact and approximate solutions
are equal to the NCSDs of two successive approximations. So we do
not need to existence  the
 exact solution and instead of applying the termination criterion (\ref{cond1}) we apply the stopping condition (\ref{cond2}).

\section{Conclusion}
This paper develops further the mathematical methodology employed in
\cite{0025}for the first multi-vector energy analysis for the
interconnected energy systems of Great Britain and Ireland. The
efficient numerical method is designed and applied to the load
leveling problem. Its efficiency is demonstrated on the real
retrospective data of the Single Electricity Market (SEM) of the
Island of Ireland. It is to be noted that employed evolutionary
dynamical models are able to take into account (in terms of Volterra
jump discontinuous kernels) both the time-dependent efficiency and
the availability of generation/storage of each energy storage
technology in the power system. The considered SEM is a power system
with high wind power penetration and much unpredictability due to
the inherent variability of wind. Following \cite{0025}, the problem
of efficient charge/discharge strategies is reduced to solving
integral equations and their systems. The novel numerical methods
using stochastic arithmetic and Taylor-collocation are proposed for
such equations to find the available storage dispatch schedules. The
convergence theorem is proved for the designed method of IE
solution. The proposed method is applied to real data demonstrating
its effectiveness.





\bigskip

\textbf{Samad Noeiaghdam}, PhD, Associate professor, Baikal School of BRICS, Irkutsk National Research
Technical University, Lermontov st. 83, Irkutsk, 664074, Russian Federation;
South Ural State University, Lenin prospect 76, Chelyabinsk, 454080, Russian Federation

\email{samadnoeiaghdam@gmail.com}

\textbf{Denis Sidorov}, Doctor of Sciences (Physics and Mathematics), Pro-fessor, Melentiev Energy Systems Institute SB RAS, 130, Lermontov st.,
Irkutsk, 664033, Russian Federation, tel.: (3952) 500-646 ext. 258; Irkutsk
National Research Technical University, Lermontov st. 83, Irkutsk, 664074,
Russian Federation; Irkutsk State University, K. Marx st. 1, Irkutsk,
664003, Russian Federation
\email{dsidorov@isem.irk.ru}

\textbf{Ildar Muftahov}, Programmer, Irkutsk Computing Center of Joint
Stock Company Russian Railways, Mayakovaskii st. 25, Irkutsk, 664005,
Russian Federation; Melentiev Energy Systems Institute SB RAS, 130,
Lermontov st., Irkutsk, 664033, Russian Federation
\email{ildar\_sm@mail.ru}

\textbf{Aleksei Zhukov}, Junior research fellow, Institute of Solar-Terrestrial Physics SB RAS, 126a, Lermontov st., Irkutsk, 664033, Russian Federation;
Melentiev Energy Systems Institute SB RAS, 130, Lermontov st., Irkutsk, 664033, Russian Federation
\email{zhukovalex13@gmail.com}


\naze{Контроль точности метода коллокаций Тейлора для задачи выравнивания нагрузки}

\avtore{Нойягдам Самад}
\inst{\\Байкальская школа БРИКС, Иркутский Национальный Исследовательский Технический Университет, г.~Иркутск, Российская Федерация\\
Южно-Уральский Государственный Университет, г.~Челябинск, Российская Федерация}

\avtore{Сидоров Денис Николаевич}
\inst{\\Институт систем энергетики им. Л.А.Мелентьева СО РАН, г.~Иркутск,~Российская Федерация\\
Иркутский национальный исследовательский технический университет, г.~Иркутск,~Российская Федерация\\
Иркутский государственный университет, г.~Иркутск,~Российская Федерация}

\avtore{Муфтахов Ильдар Ринатович}
\inst{\\Институт систем энергетики им. Л.А.Мелентьева СО РАН, г.~Иркутск, Российская Федерация\\
Иркутский информационно-вычислительный центр ОАО <<РЖД>>, \\г.~Иркутск, Российская Федерация}

\avtore{Жуков Алексей Витальевич}
\inst{\\Институт солнечно-земной физики СО РАН, г.~Иркутск, Российская Федерация\\
Институт систем энергетики им. Л.А.Мелентьева СО РАН, г.~Иркутск, Российская Федерация}

\begin{abstracte}
Высокая степень проникновения возобновляемых источников энергии в сочетании с децентрализацией транспортных и тепловых нагрузок в будущих энергосистемах приведет к еще более сложному решению проблемы энергозатрат с учетом планирования аккумулирования энергии для эффективного выравнивания нагрузки.
В данной статье рассматривается адаптивный подход к задаче выравнивания нагрузки с использованием интегральных динамических моделей Вольтерра.
Задача формулируется как решение интегрального уравнения Вольтерра первого рода, которое решается с помощью численного метода коллокаций Тейлора, имеющего точность второго порядка и обладающего свойствами саморегуляции, что связано с доверительными уровнями системного спроса.
Также применяется метод CESTAC для нахождения оптимальной аппроксимации, оптимальной погрешности и оптимального шага метода коллокаций.
Данный адаптивный подход подходит для оптимизации накопления энергии в режиме реального времени.
Эффективность предлагаемой методики продемонстрирована на едином рынке электроэнергии острова Ирландия.
\end{abstracte}

\keywordse{задача выравнивания нагрузки; метод коллокаций Тейлора; стохастическая арифметика; метод CESTAC} 

\selectlanguage{russian}

\textbf{Нойягдам Самад}, PhD, доцент,
{Иркутский национальный исследовательский технический университет},
664074, г. Иркутск, ул. Лермонтова, 83, Российская Федерация;
{Южно-Уральский Государственный Университет},
454080, г. Челябинск, Проспект Ленина, 76, Российская Федерация
\email{samadnoeiaghdam@gmail.com}.

\textbf{Сидоров Денис Николаевич}, доктор физико-математических наук, профессор РАН,
Институт систем энергетики им. Л.А.Мелентьева СО РАН,
664033, Иркутская область, г. Иркутск, ул. Лермонтова, д. 130, тел.: (3952)500-646 (код
258), Российская Федерация; Иркутский национальный исследовательский технический университет, 664074, г. Иркутск, ул. Лермонтова, 83, Российская Федерация; Иркутский государственный университет, 664003, Иркутск, ул. К. Маркса, 1, Российская Федерация
\email{contact.dns@gmail.com}.

\textbf{Муфтахов Ильдар Ринатович}, программист, Иркутский инфор-мационно-вычислительный центр ОАО <<РЖД>>, 664005, г. Иркутск, ул.
Маяковского, 25, Российская Федерация; Институт систем энергетики
им. Л. А. Мелентьева СО РАН, 664033, г. Иркутск, ул. Лермонтова, 130,
Российская Федерация
\email{ildar\_sm@mail.ru}.

\textbf{Жуков Алексей Витальевич}, младший научный сотрудник, Институт солнечно-земной физики СО РАН, 664033, Иркутская область, г. Иркутск, ул. Лермонтова, д. 126a, Российская Федерация; Институт систем энергетики
им. Л. А. Мелентьева СО РАН, 664033, г. Иркутск, ул. Лермонтова, 130,
Российская Федерация.

\email{zhukovalex13@gmail.com}

\end{article}

\end{document}